\documentclass[12pt]{amsart}
\usepackage{amsmath,amssymb,amsbsy,amsfonts,latexsym,amsopn,amstext,
                                               amsxtra,euscript,amscd,bm}
\usepackage{url}
\usepackage[colorlinks,linkcolor=blue,anchorcolor=blue,citecolor=blue]{hyperref}
\usepackage{color}
\usepackage{float}
\usepackage{verbatim}
\usepackage{bbm}

\hypersetup{breaklinks=true}

\begin{document}

\newtheorem{theorem}{Theorem}
\newtheorem{lemma}[theorem]{Lemma}
\newtheorem{claim}[theorem]{Claim}
\newtheorem{cor}[theorem]{Corollary}
\newtheorem{prop}[theorem]{Proposition}
\newtheorem{definition}{Definition}
\newtheorem{question}[theorem]{Question}
\newcommand{\hh}{{{\mathrm h}}}

\numberwithin{equation}{section}
\numberwithin{theorem}{section}
\numberwithin{table}{section}

\def\sssum{\mathop{\sum\!\sum\!\sum}}
\def\ssum{\mathop{\sum\ldots \sum}}
\def\iint{\mathop{\int\ldots \int}}

\def\squareforqed{\hbox{\rlap{$\sqcap$}$\sqcup$}}
\def\qed{\ifmmode\squareforqed\else{\unskip\nobreak\hfil
\penalty50\hskip1em\null\nobreak\hfil\squareforqed
\parfillskip=0pt\finalhyphendemerits=0\endgraf}\fi}

\newfont{\teneufm}{eufm10}
\newfont{\seveneufm}{eufm7}
\newfont{\fiveeufm}{eufm5}
%
%
\newfam\eufmfam
     \textfont\eufmfam=\teneufm
\scriptfont\eufmfam=\seveneufm
     \scriptscriptfont\eufmfam=\fiveeufm
%
%
\def\frak#1{{\fam\eufmfam\relax#1}}

\newcommand{\bflambda}{{\boldsymbol{\lambda}}}
\newcommand{\bfmu}{{\boldsymbol{\mu}}}
\newcommand{\bfxi}{{\boldsymbol{\xi}}}
\newcommand{\bfrho}{{\boldsymbol{\rho}}}

\def\fK{\mathfrak K}
\def\fT{\mathfrak{T}}

\def\fA{{\mathfrak A}}
\def\fB{{\mathfrak B}}
\def\fC{{\mathfrak C}}

\def \balpha{\bm{\alpha}}
\def \bbeta{\bm{\beta}}
\def \bgamma{\bm{\gamma}}
\def \blambda{\bm{\lambda}}
\def \bchi{\bm{\chi}}
\def \bphi{\bm{\varphi}}
\def \bpsi{\bm{\psi}}

\def\eqref#1{(\ref{#1})}

\def\vec#1{\mathbf{#1}}


\def\cA{{\mathcal A}}
\def\cB{{\mathcal B}}
\def\cC{{\mathcal C}}
\def\cD{{\mathcal D}}
\def\cE{{\mathcal E}}
\def\cF{{\mathcal F}}
\def\cG{{\mathcal G}}
\def\cH{{\mathcal H}}
\def\cI{{\mathcal I}}
\def\cJ{{\mathcal J}}
\def\cK{{\mathcal K}}
\def\cL{{\mathcal L}}
\def\cM{{\mathcal M}}
\def\cN{{\mathcal N}}
\def\cO{{\mathcal O}}
\def\cP{{\mathcal P}}
\def\cQ{{\mathcal Q}}
\def\cR{{\mathcal R}}
\def\cS{{\mathcal S}}
\def\cT{{\mathcal T}}
\def\cU{{\mathcal U}}
\def\cV{{\mathcal V}}
\def\cW{{\mathcal W}}
\def\cX{{\mathcal X}}
\def\cY{{\mathcal Y}}
\def\cZ{{\mathcal Z}}
\newcommand{\rmod}[1]{\: \mbox{mod} \: #1}

\def\cg{{\mathcal g}}

\def\vr{\mathbf r}

\def\e{{\mathbf{\,e}}}
\def\ep{{\mathbf{\,e}}_p}
\def\er{{\mathbf{\,e}}_r}

\def\Tr{{\mathrm{Tr}}}
\def\Nm{{\mathrm{Nm}}}

\def\wcI{\widetilde{\cI}}
\def\wcJ{\widetilde{\cJ}}
\def\lcm{{\mathrm{lcm}}}

\def\({\left(}
\def\){\right)}
\def\fl#1{\left\lfloor#1\right\rfloor}
\def\rf#1{\left\lceil#1\right\rceil}
\def\rnorm#1{\langle#1\rangle_r}

\def\mand{\qquad \mbox{and} \qquad}

\newcommand{\commI}[1]{\marginpar{%
\begin{color}{red}
\vskip-\baselineskip 
\raggedright\footnotesize
\itshape\hrule \smallskip I: #1\par\smallskip\hrule\end{color}}}

\newcommand{\commK}[1]{\marginpar{%
\begin{color}{blue}
\vskip-\baselineskip 
\raggedright\footnotesize
\itshape\hrule \smallskip K: #1\par\smallskip\hrule\end{color}}}




\hyphenation{re-pub-lished}

\mathsurround=1pt

\def\bfdefault{b}
\overfullrule=5pt

\def \F{{\mathbb F}}
\def \K{{\mathbb K}}
\def \Z{{\mathbb Z}}
\def \Q{{\mathbb Q}}
\def \R{{\mathbb R}}
\def \C{{\mathbb C}}
\def\Fp{\F_p}
\def \fp{\Fp^*}

\def \U{{\mathbf U}}
\def \V{{\mathbf V}}

\def\Kmn{\cK_p(m,n)}
\def\psmn{\psi_p(m,n)}
\def\SpAB{\cS_{a,p}(\cA,\cB;\cI,\cJ)}
\def\SrAB{\cS_{a,r}(\cA,\cB;\cI,\cJ)}
\def\SpkAB{\cS_{a,p^k}(\cA,\cB;\cI,\cJ)}

\def\SpA{\cS_{a,p}(\cA;\cI,\cJ)}
\def\SrA{\cS_{a,r}(\cA;\cI,\cJ)}
\def\SpkA{\cS_{a,p^k}(\cA,\cI,\cJ)}

\def\Sp{\cS_{a,p}(\cI,\cJ)}
\def\Sr{\cS_{a,r}(\cI,\cJ)}
\def\Spk{\cS_{a,p^k}(\cI,\cJ)}

\def\RpIJ{R_{a,p}(\cI,\cJ)}
\def\RrIJ{R_{a,r}(\cI,\cJ)}
\def\RrIJz{R_{a,r}(\cI,\cJ_0)}
\def\RrIJj{R_{a,r}(\cI,\cJ_j)}
\def\RpkIJ{R_{a,p^k}(\cI,\cJ)}

\def\TpIJ{T_{a,p}(\cI,\cJ)}
\def\TrIJ{T_{a,r}(\cI,\cJ)}
\def\TrIJz{T_{a,r}(\cI,\cJ_0)}
\def\TrIJj{T_{a,r}(\cI,\cJ_j)}
\def\TpkIJ{T_{a,p^k}(\cI,\cJ)}
\def \xbar{\overline x_p}

\title[Bounds for triple Exp. Sum with Exp. and Linear Function]{Bounds for triple exponential sums with mixed exponential and linear terms}

\author[Kam Hung Yau] {Kam Hung Yau}

 

\address{Department of Pure Mathematics, University of New South Wales,
Sydney, NSW 2052, Australia}
\email{kamhung.yau@unsw.edu.au}

\begin{abstract} 
We establish bounds of triple exponential sums with mixed exponential and linear terms. The method we use is by Shparlinski together with a bound of additive energy from Roche-Newton, Rudnev and Shkredov. 
\end{abstract}


\maketitle

\section{Introduction}

Particular bounds of exponential sums were first studied in Number Theory as they produce arithmetic information about certain Diophantine problems. For example, by obtaining estimates of exponential sums over primes, Vinogradov~\cite{V} was able to establish every sufficiently large odd integer can be written as a sum of three primes. Now the study of bounds for exponential sums are both for mathematical and arithmetic interest.

Let $g$ be an arbitrary integer with $\gcd(g, p) =1$. We denote $T$ to be the multiplicative order of $g$ modulo $p$. Given two intervals of consecutive integers
$$
\cI = \{ K+1, \ldots, K+M \}, \quad \cJ = \{ L+1, \ldots, L+N \}
$$
and
$$
\cK = \{ 1, \ldots, H \}
$$
with integers $H,K,L,M,N$ such that $0 < M \le p $, $0< N \le T$, $0 < H <T$ and a complex sequence $\cA = (\alpha_{m})_{m \in \cI}$, we define the following exponential sum
$$
\cS_{a,T,p}(\cA; \cI, \cJ, \cK ) = \sum_{m \in \cI} \sum_{n \in \cJ} \sum_{x \in \cK}  \alpha_m   e_p(am g^x) e_T(nx)
$$
for integers $a \in \mathbb{Z}$ with $\gcd(a,p)=1$ where $e_{h}(x)=e(2 \pi i x/h)$. In particular, when $\cI = \mathbb{Z}_{p}$, we define
$$
\cS_{a,T,p}(\cA;  \cJ, \cK ) = \cS_{a,T,p}(\cA; \cI, \cJ, \cK ).
$$

Similar double exponential sums has already been considered. In particular, sums of the form
$$
S(\cA, \cB; \cI , \cJ) = \sum_{m \in  \cI } \sum_{n \in \cJ} \alpha_{m} \beta_{n} e_{p}(am g^{n})
$$
has been considered in the work by Shparlinski \& Yau~\cite{SY}. For the case when $g$ is not necessary a primitive root of $p$, bounds has been established under the condition $ \cI = \{ 1 \}$ and $\alpha_{m} =\beta_{n} =1$ by Kerr~\cite{K} but the same method imployed there also works for the general $\cI$ as the bound depend only on the norm. Similar sums for multiplicative character has also been studied in~\cite{SY2}. We refer the reader to~\cite{KS} for a broader overview of this subject.

In this paper we establish bounds for $\cS_{a,T,p}(\cA; \cI, \cJ, \cK )$ when $\cI = \mathbb{Z}_{p}$, it is clear the same method also works for general $\cI$.

Our approach follows from Shparlinski as in the proof of~\cite{Shp}[Theorem 2.1]. In particular, after applying the triangle and H\"older inequality to $\cS_{a,T,p}(\cA; \cI, \cJ, \cK )$, we obtain a mean fourth-moment of an exponential sum. By opening and changing the order of summation and appealing to the orthogonality of the exponential function, we can bound the sum by the number of solutions to a particular congruence (see Lemma~\ref{additive energy}).

\section{Main Result}

The statement $A \ll B$ and $A = O(B)$ are both equivalent to the inequality $|A| \le c B$ for some positive absolute constant $c$. We define for any real number $\sigma >0$,
$$
\lVert \cA \rVert_{\sigma} = \Big (\sum_{m \in \cI} | \alpha_{m}|^{\sigma} \Big )^{1/\sigma}.
$$

We state below a bound for $\cS_{a,T,p}(\cA; \cK, \cJ)$.

\begin{theorem} \label{S-bound}
For any prime $p$, we have
\begin{equation*}
\begin{split} 
S_{a,T,p}&(\cA; \cJ, \cK ) \ll  \lVert \cA \rVert_{1}^{1/2} \lVert \cA \rVert_{2}^{1/2}p^{1/4} N^{3/8} T^{5/8}.
\end{split}
\end{equation*}
\end{theorem}

Using the same technique as in~\cite[Lemma 3.14]{Shp2} and the bound~\cite[Corollary 19]{R-NRS}, we obtain the trivial bound
\begin{equation}
\begin{split} \label{trivial-bound1}
S_{a,T,p}(\cA; \cJ, \cK ) & \ll  \lVert A \rVert_{1} N \min \{ p^{1/8} H^{5/8}, p^{1/4} H^{3/8} \}.
\end{split}
\end{equation}
Assuming $| \alpha_{m}| \le 1$ we have $\lVert \cA \rVert_{1} \ll M$ and $\lVert \cA \rVert_{2} \ll M^{1/2}$. We see that Theorem~\ref{S-bound} provides a stronger bound 
$$
S_{a,T,p}(\cA; \cJ, \cK ) \ll M^{3/4}p^{1/4} N^{3/8}  T^{5/8}
$$
than~(\ref{trivial-bound1}) which becomes
$$
S_{a,T,p}(\cA; \cJ, \cK ) \ll M N \min \{ p^{1/8} H^{5/8}, p^{1/4} H^{3/8} \}
$$
when
$$
pT^{5} < M^{2}N^{5}H^{5} \mand  T^{5} < M^{2}N^{5}H^{3}.
$$

\section{Preparation}
For an integer $u$, we define
$$
\langle u \rangle_{r} = \min_{k \in \mathbb{Z}} | u - kr|
$$
as the distance to the nearest integral multiple of $r$.

We recall a well-known bound from~\cite[Bound (8.6)]{IK}.
\begin{lemma} \label{linearbound}
For an integers $u$, $W$ and $Z \ge 1$, we have
$$
\sum_{n=W+1}^{W+Z} e_{r}(nu) \ll \min \left \{ Z, \frac{r}{ \langle  u \rangle_{r}} \right \}.
$$
\end{lemma}

We recall that $T$ is the multiplicative order of $g$ modulo $p$.
For any positive integer $K \le T$, we define the \textit{additive energy} $E_{p}(K)$ as the number of solutions to the congruence
\begin{equation} \label{addenergy}
g^{x_{1}} + g^{x_{2}} \equiv g^{x_{3}} + g^{x_{4}} \pmod{p}
\end{equation}
where
$$
 \quad (x_{1}, x_{2}, x_{3}, x_{4}) \in \{1, \ldots, K \}^{4}.
$$
Our approach to bounding $\cS_{a,T,p}(\cA; \cI, \cJ, \cK )$ is to reduce the problem to estimating $E_{p}(K)$.

Note that $(v_{1}, v_{2}, v_{1}, v_{2}) \in \{1, \ldots, K \}^{4}$ is always a solution to~(\ref{addenergy}), hence we have the trivial lower bound $K^{2} \le E_{p}(K)$. If $(v_{1}, v_{2}, v_{3}, v_{4}) \in \{1, \ldots, K \}^{4}$ is a solution to~(\ref{addenergy}) then $v_{4}$ is depended on $v_{1}, v_{2}, v_{3}$ and we obtain the trivial upper bound $E_{p}(K) \le K^3$. In particular, $E_{p}(K)$ is an increasing function of $K$.

Set $A, B, C= \{g, \ldots, g^{K} \}$ then we have the trivial bound $|A|\le K$,  $|BC| \le 2K$. Appealing to~\cite[Theorem 6]{R-NRS}, we derive a non-trivial estimate on $E_{p}(K)$. 

\begin{lemma} \label{additive energy}
For any positive integer $1 \le K \le T$, we have
$$
E_{p}(K) \ll K^{5/2}.
$$
\end{lemma}

\section{Proof of Theorem~\ref{S-bound}}
We proceed similarly to the proof of~\cite[Theorem 2.1]{Shp}.
Rearranging then applying Lemma~\ref{linearbound}, we have 
\begin{align*}
\cS_{a,T,p}(\cA; \cJ, \cK ) 
& = \sum_{x=1}^{H} \sum_{m =0}^{p-1} \alpha_{m} e_{p}(amg^{x}) \sum_{n=L+1}^{L+N} e_{T}(nx)  \\
& =   \sum_{x=1}^{H} \sum_{m = 0}^{p-1} \alpha_{m} e_{p}(amg^{x}) \varphi_{x}
\end{align*}
where
$$
| \varphi_{x}  | \le \min \left (  N, \frac{T}{ \langle x \rangle_{T}} \right ).
$$
Define $I = \lceil \log N \rceil$ and define the sets
$$
\cL_{0} = \{ x \in \mathbb{Z}: 0 <  x \le  T/N \}
$$
and
$$
\cL_{i} = \{ x \in \mathbb{Z}: \min \{T ,e^{i}T/N \}  \ge  x >  e^{i-1}T/N \}
$$
for $i=1, \ldots, I$. Therefore, we obtain
$$
\cS_{a,T,p}(\cA; \cJ , \cK)  \ll \sum_{i=0}^{I} |S_{i}|
$$
where
$$
S_{i} = \sum_{x \in \cL_{i}} \sum_{m=0}^{p-1} \alpha_{m} e_{p}(amg^{x}) \varphi_{x} 
$$
for $i=0, \ldots, I.$

Applying the triangle and H\"older inequality, we obtain
\begin{align}   \label{Si}
|S_{i}|  &\le  \sum_{m=0}^{p-1} |\alpha_{m}|^{1/2} |\alpha_{m}^{2}|^{1/4} \Big | \sum_{x \in \cL_{i}} \alpha_{m} e_{p}(amg^{x}) \varphi_{x}  \Big | \nonumber \\
& \le  \Big ( \sum_{m=0}^{p-1} |\alpha_{m}|  \Big )^{1/2}  \Big ( \sum_{m=0}^{p-1} | \alpha_{m}|^{2}  \Big )^{1/4}  \Big ( \sum_{m=0}^{p-1}  \Big |  \sum_{x \in \cL_{i}}  e_{p}(amg^{x}) \varphi_{x}  \Big |^{4} \Big )^{1/4}  \\
& = \lVert \cA \rVert_{1}^{1/2} \lVert \cA \rVert_{2}^{1/2} \Big (  \sum_{m=0}^{p-1}  \Big |  \sum_{x \in \cL_{i}}  e_{p}(amg^{x}) \varphi_{x}  \Big  |^{4} \Big )^{1/4}   \nonumber
\end{align}
which is valid for all $i = 0, \ldots , I$. Opening the summation and changing the order of summation, we obtain
\begin{align*}
\sum_{m=0}^{p-1} &  \Big |  \sum_{x \in \cL_{i}}  e_{p}(amg^{x}) \varphi_{x}  \Big |^{4} \\
& = \sum_{m=0}^{p-1} \underset{x_{1}, \ldots,x_{4} \in \cL_{i}}{ \sum \cdots \sum} \varphi_{x_{1}} \varphi_{x_{2}} \overline{ \varphi_{x_{3}} \varphi_{x_{4 }} } e_{p}(am(g^{x_{1}} + g^{x_{2}} - g^{x_{3}} -g^{x_{4}}    )) \\
& = \underset{x_{1}, \ldots,x_{4} \in \cL_{i}}{ \sum \cdots \sum} \varphi_{x_{1}} \varphi_{x_{2}} \overline{ \varphi_{x_{3}} \varphi_{x_{4 }} }   \sum_{m=0}^{p-1}  e_{p}(am(g^{x_{1}} + g^{x_{2}} - g^{x_{3}} -g^{x_{4}}    )).
\end{align*}
Since for all $x \in \cL_{i}$, we have the bound $\varphi_{x} \ll e^{-i}N$, hence we get
\begin{align*}
\sum_{m=0}^{p-1} &   \Big  |  \sum_{x \in \cL_{i}}  e_{p}(amg^{x}) \varphi_{x}  \Big  |^{4} \\
& \le   \underset{x_{1}, \ldots,x_{4} \in \cL_{i}}{ \sum \cdots \sum} |\varphi_{x_{1}} \varphi_{x_{2}} \overline{ \varphi_{x_{3}} \varphi_{x_{4 }} } |\sum_{m=0}^{p-1}  e_{p}(am(g^{x_{1}} + g^{x_{2}} - g^{x_{3}} -g^{x_{4}}    )) \\
& \ll  e^{-4i}N^{4} \underset{x_{1}, \ldots,x_{4} \in \cL_{i}}{ \sum \cdots \sum} \sum_{m=0}^{p-1}  e_{p}(am(g^{x_{1}} + g^{x_{2}} - g^{x_{3}} -g^{x_{4}}    )). 
\end{align*}
By appealing to the orthogonality of exponential function, we obtain 
$$
\sum_{m=0}^{p-1}    \Big  |  \sum_{x \in \cL_{i}}  e_{p}(amg^{x}) \varphi_{x}  \Big  |^{4} \ll p  e^{-4i}N^{4}  E_{p}( \lfloor e^{i} T/N \rfloor).
$$
Therefore by Lemma~\ref{additive energy}, we obtain
\begin{align*}
\sum_{m=0}^{p-1}   \left |  \sum_{x \in \cL_{i}} e_{p}(amg^{x}) \varphi_{x}  \right |^{4} & \ll   p  e^{-4i}N^{4} (e^{i} T/N)^{5/2} \\
& \ll p e^{-3/2i} N^{3/2}T^{5/2}.
\end{align*}
Substituting this bound into~(\ref{Si}), we obtain
$$
|S_{i}| \ll \lVert \cA \rVert_{1}^{1/2} \lVert \cA \rVert_{2}^{1/2}p^{1/4}  e^{-3i/8}N^{3/8} T^{5/8}. 
$$
Finally, we have
$$
\sum_{i=0}^{I} |S_{i}| \ll \lVert \cA \rVert_{1}^{1/2} \lVert \cA \rVert_{2}^{1/2}p^{1/4} N^{3/8} T^{5/8}
$$
and the result follows immediately.

\section*{Acknowledgement}
The author would like to thank the referee for many helpful comments and also I. E. Shparlinski for the problem, helpful comments and proof-reading of this paper. During the preparation of this paper, the author was supported by an Australian
Government Research Training Program (RTP) Scholarship.

\end{document}